\documentclass{amsart}
\usepackage{amsmath}
\usepackage{amssymb}
\usepackage{amsfonts}

\newtheorem{theorem}{Theorem}
\theoremstyle{plain}

\newtheorem{remark}{Remark}

\numberwithin{equation}{section}
\input{tcilatex}

\begin{document}
\title[Sampling Theorems for Sturm Liouville Problem]{Sampling Theorems for
Sturm Liouville Problem with Moving Discontinuity Points}
\author{Fatma H\i ra}
\address[F. H\i ra \ and N. Alt\i n\i \c{s}\i k]{The University of Ondokuz
May\i s, Science and Arts Faculty, Department of Mathematics,Samsun, Turkey\\
}
\email[Corresponding author: F. H\i ra]{fatma.hira@omu.edu.tr}
\urladdr{}
\thanks{}
\author{Nihat Alt\i n\i \c{s}\i k}
\curraddr{ }
\email[N.~Alt\i n\i \c{s}\i k]{anihat@omu.edu.tr}
\address{ }
\subjclass[2000]{ 34B24, 34B27; 94A20}
\keywords{Symmetric and Moving Discontinuities; Discontinuous
Sturm-Liouville Problem; Green's Function; Sampling Theory}

\begin{abstract}
In this paper, we investigate the sampling analysis for a new
Sturm-Liouville problem with symmetrically located discontinuities which are
defined to depending on a neighborhood of a midpoint of the interval. Also
the problem has transmission conditions at these points of discontinuity and
includes an eigenparameter in a boundary condition. We establish briefly the
needed relations for the derivations of the sampling theorems and construct
Green's function for the problem. Then we derive sampling representations
for transforms whose kernels are either solutions or Green's functions.
\end{abstract}

\maketitle

\section{Introduction}

Throughout this study we consider Sturm Liouville problem:%
\begin{equation}
\tau \left( u\right) :=-u^{\shortparallel }\left( x\right) +q\left( x\right)
u\left( x\right) =\lambda u\left( x\right) ,~~\ \ \ \ x\in I,  \tag{1.1}
\label{1.1}
\end{equation}%
with an eigenparameter dependent on a boundary condition;%
\begin{equation}
B_{a}\left( u\right) :=\beta _{1}u\left( a\right) +\beta _{2}u^{\prime
}\left( a\right) =0,  \tag{1.2}  \label{1.2}
\end{equation}%
\begin{equation}
B_{b}\left( u\right) :=\lambda \left( \alpha _{1}^{\prime }u\left( b\right)
-\alpha _{2}^{\prime }u^{\prime }\left( b\right) \right) +\left( \alpha
_{1}u\left( b\right) -\alpha _{2}u^{\prime }\left( b\right) \right) =0, 
\tag{1.3}  \label{1.3}
\end{equation}%
and transmission conditions at two points of discontinuity depending on a
neighborhood of $\ \theta ,$ that are $\theta _{-\varepsilon }$ and $\theta
_{+\varepsilon }$;%
\begin{equation}
T_{-\varepsilon }\left( u\right) :=u\left( \theta _{-\varepsilon }-\right)
-\delta u\left( \theta _{-\varepsilon }+\right) =0,  \tag{1.4}  \label{1.4}
\end{equation}%
\begin{equation}
T_{-\varepsilon }^{\prime }\left( u\right) :=u^{\prime }\left( \theta
_{-\varepsilon }-\right) -\delta u^{\prime }\left( \theta _{-\varepsilon
}+\right) =0,  \tag{1.5}  \label{1.5}
\end{equation}%
\begin{equation}
T_{+\varepsilon }\left( u\right) :=\delta u\left( \theta _{+\varepsilon
}-\right) -\gamma u\left( \theta _{+\varepsilon }+\right) =0,  \tag{1.6}
\label{1.6}
\end{equation}%
\begin{equation}
T_{+\varepsilon }^{\prime }\left( u\right) :=\delta u^{\prime }\left( \theta
_{+\varepsilon }-\right) -\gamma u^{\prime }\left( \theta _{+\varepsilon
}+\right) =0,  \tag{1.7}  \label{1.7}
\end{equation}%
where $I=\left[ a,\theta _{-\varepsilon }\right) \cup \left( \theta
_{-\varepsilon },\theta _{+\varepsilon }\right) \cup \left( \theta
_{+\varepsilon },b\right] ;$ $\lambda $ is a complex spectral parameter, $%
q\left( x\right) $ is a given real valued function which is continuous in $%
\left[ a,\theta _{-\varepsilon }\right) ,~\left( \theta _{-\varepsilon
},\theta _{+\varepsilon }\right) ,~$and $\left( \theta _{+\varepsilon },b%
\right] $ and has a finite limit $q\left( \theta _{-\varepsilon }\pm \right) 
$ and $q\left( \theta _{+\varepsilon }\pm \right) $; $\beta _{i},\alpha
_{i},\alpha _{i}^{\prime },\delta ,\gamma \in 
\mathbb{R}
$ $\left( i=1,2\right) $; $\left\vert \beta _{1}\right\vert +\left\vert
\beta _{2}\right\vert \neq 0,$ $\delta \neq 0,$ $\gamma \neq 0;$ $\theta
:=\left( a+b\right) /2,$ $\theta _{\pm \varepsilon }\pm :=\left( \theta \pm
\varepsilon \right) \pm 0,$ $0<\varepsilon <\left( b-a\right) /2$ and%
\begin{equation}
\rho :=\left( \alpha _{1}^{\prime }\alpha _{2}-\alpha _{1}\alpha
_{2}^{\prime }\right) >0.  \tag{1.8}  \label{1.8}
\end{equation}

In the literature, the Whittaker-Kotel'nikov-Shannon (WKS) sampling theorem
and generalization of the WKS sampling theorem (see$\left[ 1-3\right] )$ has
been investigated extensively (see $\left[ 4-8\right] ).$ Sampling theorems
associated with Sturm-Liouville problems were investigated in $\left[ 9-13%
\right] .$ Also, $\left[ 14-17\right] $ and $\left[ 18-21\right] $ are the
examples works in direction of sampling analysis associated with continuous
and discontinuous eigenproblems, respectively. In $\left[ 20\right] $ the
author investigated the sampling analysis associated with discontinuous
Sturm Liouville problems which has transmission conditions at the point of
discontinuity and contains an eigenparameter in two boundary canditions. In
the present study, we introduce a new Sturm Liouville problem which has
symmetrically located discontinuities which are defined to depending on a
neighborhood of a midpoint of the interval. $\varepsilon $ is a parameter
controling the change of neighborhood process (it can be called tuning
parameter) and by using the change of this $\varepsilon $ parameter it's
possible to determine points of discontinuity. That is, two points of
discontinuity can be determined in interval $\left[ a,b\right] $ for each $%
\varepsilon $ value in interval $0<\varepsilon <\left( b-a\right) /2.$ This
is the difference between the problem $\left( 1.1\right) -\left( 1.7\right) $
and Sturm Liouville eigenvalue problem studied extensively in the literature
(see $\left[ 22-28\right] ).$ The main result is that points of
discontinuity can be determined and moved by changing $\varepsilon $
parameter. A similar problem in more detail was presented in $\left[ 28%
\right] .$ This study is the first to investigate that sampling analysis
associated with eigenproblems with moving discontinuity points . To derive
sampling theorems for the problem $\left( 1.1\right) -\left( 1.7\right) ,$
we establish briefly the spectral properties and construct Green's function
of the problem $\left( 1.1\right) -\left( 1.7\right) $. Then we derive two
sampling theorems using solutions and Green's function, respectively.

\section{An Operator Formulation and Asymptotic Formulas}

To formulate a theoretic approach to the problem (1.1)-(1.7) we define the
Hilbert space $H=L_{2}\left( a,b\right) \oplus 
\mathbb{C}
$ with an inner product%
\begin{equation}
\left\langle \text{f}\left( .\right) ,\text{g}\left( .\right) \right\rangle
_{H}:=\underset{a}{\overset{\theta _{-\varepsilon }}{\int }}f\left( x\right) 
\overline{g}\left( x\right) dx+\delta ^{2}\underset{\theta _{-\varepsilon }}{%
\overset{\theta _{+\varepsilon }}{\int }}f\left( x\right) \overline{g}\left(
x\right) dx+\gamma ^{2}\underset{\theta _{+\varepsilon }}{\overset{b}{\int }}%
f\left( x\right) \overline{g}\left( x\right) dx+\frac{\gamma ^{2}}{\rho }h%
\overline{k},  \tag{2.1}  \label{2.1}
\end{equation}%
where f$\left( x\right) =\left( 
\begin{array}{c}
f\left( x\right) \\ 
h%
\end{array}%
\right) ,$ g$\left( x\right) =\left( 
\begin{array}{c}
g\left( x\right) \\ 
k%
\end{array}%
\right) \in H,$ $f\left( .\right) ,g\left( .\right) \in L_{2}\left(
a,b\right) ,$ $h,k\in 
\mathbb{C}
.$ For convenience we put%
\begin{equation}
R\left( u\right) :=\alpha _{1}u\left( b\right) -\alpha _{2}u^{\prime }\left(
b\right) ,\text{ \ \ }R^{\prime }\left( u\right) :=\alpha _{1}^{\prime
}u\left( b\right) -\alpha _{2}^{\prime }u^{\prime }\left( b\right) . 
\tag{2.2}  \label{2.2}
\end{equation}

Let $D\left( A\right) \subseteq H$ be the set of all f$\left( x\right)
=\left( 
\begin{array}{c}
f\left( x\right) \\ 
h%
\end{array}%
\right) \in H$ such that $f$ and $f^{\prime }$ are absolutely continuous on $%
\left[ a,b\right] $ and $\tau \left( f\right) \in L_{2}\left( a,b\right) ,$ $%
h=R^{\prime }\left( f\right) ,$ $B_{a}\left( f\right) =0,$ $T_{\pm
\varepsilon }\left( f\right) =T_{\pm \varepsilon }^{\prime }\left( f\right)
=0.$ Define the operator $A:D\left( A\right) \rightarrow H$ by%
\begin{equation}
A\left( 
\begin{array}{c}
f\left( x\right) \\ 
R^{\prime }\left( f\right)%
\end{array}%
\right) =\left( 
\begin{array}{c}
\tau \left( f\right) \\ 
-R\left( f\right)%
\end{array}%
\right) ,\text{ \ \ }\left( 
\begin{array}{c}
f\left( x\right) \\ 
R^{\prime }\left( f\right)%
\end{array}%
\right) \in D\left( A\right) .  \tag{2.3}  \label{2.3}
\end{equation}

The operator $A:D\left( A\right) \rightarrow H$ is equivalent to the
eigenvalue problem (1.1)-(1.7) in the sense that the eigenvalues of $A$ are
exactly those of the problem (1.1)-(1.7).

We can prove in a manner similar to that of $\left[ 23,25,26,28\right] $
that $A$ is symmetric in $H,$ all eigenvalues of the problem are real.

Let $\phi _{\lambda }\left( .\right) $ and $\chi _{\lambda }\left( .\right) $
be two solutions of (1.1) as%
\begin{equation}
\phi _{\lambda }\left( x\right) =\left\{ 
\begin{array}{l}
\phi _{-\varepsilon ,\lambda }\left( x\right) ,\text{ \ \ }x\in \left[
a,\theta _{-\varepsilon }\right) , \\ 
\phi _{\varepsilon ,\lambda }\left( x\right) ,\text{ \ }x\in \left( \theta
_{-\varepsilon },\theta _{+\varepsilon }\right) , \\ 
\phi _{+\varepsilon ,\lambda }\left( x\right) ,\text{ \ \ }x\in \left(
\theta _{+\varepsilon },b\right] ,%
\end{array}%
\right. \text{ \ \ \ }\chi _{\lambda }\left( x\right) =\left\{ 
\begin{array}{l}
\chi _{-\varepsilon ,\lambda }\left( x\right) ,\text{ \ \ }x\in \left[
a,\theta _{-\varepsilon }\right) , \\ 
\chi _{\varepsilon ,\lambda }\left( x\right) ,\text{ \ }x\in \left( \theta
_{-\varepsilon },\theta _{+\varepsilon }\right) , \\ 
\chi _{+\varepsilon ,\lambda }\left( x\right) ,\text{ \ \ }x\in \left(
\theta _{+\varepsilon },b\right] ,%
\end{array}%
\right.  \tag{2.4}  \label{2.4}
\end{equation}%
satisfying the following conditions, respectively;%
\begin{equation}
\phi _{-\varepsilon ,\lambda }\left( a\right) =\beta _{2},\text{\ \ \ \ \ }%
\phi _{-\varepsilon ,\lambda }^{\prime }\left( a\right) =-\beta _{1}, 
\tag{2.5}  \label{2.5}
\end{equation}%
\begin{equation}
\phi _{\varepsilon ,\lambda }\left( \theta _{-\varepsilon }\right) =\delta
^{-1}\phi _{-\varepsilon ,\lambda }\left( \theta _{-\varepsilon }-\right) ,%
\text{ \ \ \ \ \ }\phi _{\varepsilon ,\lambda }^{\prime }\left( \theta
_{-\varepsilon }\right) =\delta ^{-1}\phi _{-\varepsilon ,\lambda }^{\prime
}\left( \theta _{-\varepsilon }-\right) ,  \tag{2.6}  \label{2.6}
\end{equation}%
\begin{equation}
\phi _{+\varepsilon ,\lambda }\left( \theta _{+\varepsilon }\right) =\delta
\gamma ^{-1}\phi _{\varepsilon ,\lambda }\left( \theta _{+\varepsilon
}-\right) ,\text{ \ \ \ \ }\phi _{+\varepsilon ,\lambda }^{\prime }\left(
\theta _{+\varepsilon }\right) =\delta \gamma ^{-1}\phi _{\varepsilon
,\lambda }^{\prime }\left( \theta _{+\varepsilon }-\right) ,  \tag{2.7}
\label{2.7}
\end{equation}%
and%
\begin{equation}
\chi _{+\varepsilon ,\lambda }\left( b\right) =\lambda \alpha _{2}^{\prime
}+\alpha _{2},\text{ \ \ \ \ }\chi _{+\varepsilon ,\lambda }^{\prime }\left(
b\right) =\lambda \alpha _{1}^{\prime }+\alpha _{1}\text{,}  \tag{2.8}
\label{2.8}
\end{equation}%
\begin{equation}
\chi _{\varepsilon ,\lambda }\left( \theta _{+\varepsilon }\right) =\gamma
\delta ^{-1}\chi _{+\varepsilon ,\lambda }\left( \theta _{+\varepsilon
}+\right) ,\text{ \ \ \ }\chi _{\varepsilon ,\lambda }^{\prime }\left(
\theta _{+\varepsilon }\right) =\gamma \delta ^{-1}\chi _{+\varepsilon
,\lambda }^{\prime }\left( \theta _{+\varepsilon }+\right) ,  \tag{2.9}
\label{2.9}
\end{equation}%
\begin{equation}
\chi _{-\varepsilon ,\lambda }\left( \theta _{-\varepsilon }\right) =\delta
\chi _{\varepsilon ,\lambda }\left( \theta _{-\varepsilon }+\right) ,\text{
\ \ \ }\chi _{-\varepsilon ,\lambda }^{\prime }\left( \theta _{-\varepsilon
}\right) =\delta \chi _{\varepsilon ,\lambda }^{\prime }\left( \theta
_{-\varepsilon }+\right) .  \tag{2.10}  \label{2.10}
\end{equation}%
These functions are entire in $\lambda $ for all $x\in \left[ a,b\right] .$

Let $W\left( \phi _{\lambda },\chi _{\lambda };x\right) $ be the Wronskian
of $\phi _{\lambda }$ and $\chi _{\lambda }$ which is independent of $x,$
since the coefficient of $y^{\prime }$ in the equation $\left( 1.1\right) $
is zero. Let \ \ 
\begin{eqnarray}
\omega \left( \lambda \right) &:&=\text{ \ }W\left( \phi _{\lambda },\chi
_{\lambda };x\right) =\phi _{\lambda }\left( x\right) \chi _{\lambda
}^{\prime }\left( x\right) -\phi _{\lambda }^{\prime }\left( x\right) \chi
_{\lambda }\left( x\right) \text{ \ }  \notag \\
\text{ \ \ \ \ \ \ } &=&\omega _{-\varepsilon }\left( \lambda \right)
=\delta ^{2}\omega _{\varepsilon }\left( \lambda \right) =\gamma ^{2}\omega
_{+\varepsilon }\left( \lambda \right) .\text{ \ }  \TCItag{2.11}
\label{2.11}
\end{eqnarray}

Then $\omega \left( \lambda \right) $ is an entire function of $\lambda $
whose zeros are precisely the eigenvalues of the operator $A.$ Using
techniques similar of those established by Titchmarsh in $\left[ 22\right] ,$
see also $\left[ 25,26,28\right] $ the zeros of $\omega \left( \lambda
\right) $ are real and simple and if $\lambda _{n},n=0,1,2,...$ denote the
zeros of $\omega \left( \lambda \right) ,$ then the two component vectors%
\begin{equation}
\Phi _{n}\left( x\right) :=\left( 
\begin{array}{c}
\phi _{\lambda _{n}}\left( x\right) \\ 
R^{\prime }\left( \phi _{\lambda _{n}}\right)%
\end{array}%
\right)  \tag{2.12}  \label{2.12}
\end{equation}%
are the corresponding eigenvectors of the operator $A$ satisfying the
orthogonality relation%
\begin{equation}
\left\langle \Phi _{n}\left( .\right) ,\Phi _{m}\left( .\right)
\right\rangle _{H}=0,\text{ for }n\neq m.  \tag{2.13}  \label{2.13}
\end{equation}

Here $\left\{ \phi _{\lambda _{n}}\left( .\right) \right\} _{n=0}^{\infty }$
will be the sequence of eigenfunctions of the problem (1.1)-(1.7)
corresponding to the eigenvalues $\left\{ \lambda _{n}\right\}
_{n=0}^{\infty }$ . We denote by%
\begin{equation}
\Psi _{n}\left( x\right) :=\frac{\Phi _{n}\left( x\right) }{\left\Vert \Phi
_{n}\left( .\right) \right\Vert _{H}}=\left( 
\begin{array}{c}
\Psi _{n}\left( x\right) \\ 
R^{\prime }\left( \Psi _{n}\right)%
\end{array}%
\right) .  \tag{2.14}  \label{2.14}
\end{equation}

Let $k_{n}\neq 0$ be the real constants for which%
\begin{equation}
\chi _{\lambda _{n}}\left( x\right) =k_{n}\phi _{\lambda _{n}}\left(
x\right) ,\text{ \ \ \ }x\in I,\text{ }n=0,1,2,....  \tag{2.15}  \label{2.15}
\end{equation}

The asymptotics of the eigenvalues and eigenfunctions can be derived similar
to the classical techniques of $\left[ 23,25,26,28\right] .$ We state the
results briefly.

$\phi _{\lambda }\left( .\right) $ is the solution determined by equations
(2.5)-(2.7) above then the following integral equations hold for $k=0$ and $%
k=1:$%
\begin{eqnarray}
\dfrac{d^{k}}{dx^{k}}\phi _{-\varepsilon ,\lambda }\left( x\right) &=&\beta
_{2}\dfrac{d^{k}}{dx^{k}}\left( \cos \sqrt{\lambda }\left( x-a\right)
\right) -\frac{\beta _{1}}{\sqrt{\lambda }}\dfrac{d^{k}}{dx^{k}}(\sin \sqrt{%
\lambda }\left( x-a\right) +  \notag \\
&&\frac{1}{\sqrt{\lambda }}\underset{a}{\overset{x}{\int }}\dfrac{d^{k}}{%
dx^{k}}\left( \sin \sqrt{\lambda }\left( x-y\right) \right) q\left( y\right)
\phi _{-\epsilon ,\lambda }\left( y\right) dy,  \TCItag{2.16}  \label{2.16}
\end{eqnarray}%
\begin{eqnarray}
\dfrac{d^{k}}{dx^{k}}\phi _{\varepsilon ,\lambda }\left( x\right) &=&\delta
^{-1}\phi _{-\epsilon ,\lambda }\left( \theta _{-\epsilon }-\right) \dfrac{%
d^{k}}{dx^{k}}\left( \cos \sqrt{\lambda }\left( x-\theta _{-\epsilon
}\right) \right) -  \notag \\
&&\frac{\delta ^{-1}}{\sqrt{\lambda }}\phi _{-\epsilon ,\lambda }^{\prime
}\left( \theta _{-\epsilon }-\right) \dfrac{d^{k}}{dx^{k}}(\sin \sqrt{%
\lambda }\left( x-\theta _{-\epsilon }\right) +  \notag \\
&&\frac{1}{\sqrt{\lambda }}\underset{\theta _{-\varepsilon }}{\overset{x}{%
\int }}\dfrac{d^{k}}{dx^{k}}\left( \sin \sqrt{\lambda }\left( x-y\right)
\right) q\left( y\right) \phi _{\epsilon ,\lambda }\left( y\right) dy, 
\TCItag{2.17}  \label{2.17}
\end{eqnarray}%
\begin{eqnarray}
\dfrac{d^{k}}{dx^{k}}\phi _{+\varepsilon ,\lambda }\left( x\right) &=&\delta
\gamma ^{-1}\phi _{\varepsilon ,\lambda }\left( \theta _{+\varepsilon
}-\right) \dfrac{d^{k}}{dx^{k}}\left( \cos \sqrt{\lambda }\left( x-\theta
_{+\epsilon }\right) \right) -  \notag \\
&&\frac{\delta \gamma ^{-1}}{\sqrt{\lambda }}\phi _{\varepsilon ,\lambda
}^{\prime }\left( \theta _{+\varepsilon }-\right) \dfrac{d^{k}}{dx^{k}}(\sin 
\sqrt{\lambda }\left( x-\theta _{+\epsilon }\right) +  \notag \\
&&\frac{1}{\sqrt{\lambda }}\underset{\theta _{+\varepsilon }}{\overset{x}{%
\int }}\dfrac{d^{k}}{dx^{k}}\left( \sin \sqrt{\lambda }\left( x-y\right)
\right) q\left( y\right) \phi _{+\varepsilon ,\lambda }\left( y\right) dy, 
\TCItag{2.18}  \label{2.18}
\end{eqnarray}%
sufficiently large $\lambda $ and $\phi _{\lambda }\left( .\right) $ have
the following asymptotic representations for $\left\vert \lambda \right\vert
\rightarrow \infty ,$ which hold uniformly for $x\in I$:%
\begin{equation}
\dfrac{d^{k}}{dx^{k}}\phi _{-\varepsilon ,\lambda }\left( x\right) =\beta
_{2}\dfrac{d^{k}}{dx^{k}}\left( \cos \sqrt{\lambda }\left( x-a\right)
\right) +O\left( \left( \sqrt{\lambda }\right) ^{k-1}e^{\left\vert
t\right\vert \left( x-a\right) }\right) ,  \tag{2.19}  \label{2.19}
\end{equation}%
\begin{equation}
\dfrac{d^{k}}{dx^{k}}\phi _{\varepsilon ,\lambda }\left( x\right) =\beta
_{2}\delta ^{-1}\dfrac{d^{k}}{dx^{k}}\left( \cos \sqrt{\lambda }\left(
x-a\right) \right) +O\left( \left( \sqrt{\lambda }\right)
^{k-1}e^{\left\vert t\right\vert \left( x-a\right) }\right) ,  \tag{2.20}
\label{2.20}
\end{equation}%
\begin{equation}
\dfrac{d^{k}}{dx^{k}}\phi _{+\varepsilon ,\lambda }\left( x\right) =\beta
_{2}\gamma ^{-1}\dfrac{d^{k}}{dx^{k}}\left( \cos \sqrt{\lambda }\left(
x-a\right) \right) +O\left( \left( \sqrt{\lambda }\right)
^{k-1}e^{\left\vert t\right\vert \left( x-a\right) }\right) ,  \tag{2.21}
\label{2.21}
\end{equation}%
if $\beta _{2}\neq 0,$%
\begin{equation}
\dfrac{d^{k}}{dx^{k}}\phi _{-\varepsilon ,\lambda }\left( x\right) =-\frac{%
\beta _{1}}{\sqrt{\lambda }}\dfrac{d^{k}}{dx^{k}}\left( \sin \sqrt{\lambda }%
\left( x-a\right) \right) +O\left( \left( \sqrt{\lambda }\right)
^{k-2}e^{\left\vert t\right\vert \left( x-a\right) }\right) ,  \tag{2.22}
\label{2.22}
\end{equation}%
\begin{equation}
\dfrac{d^{k}}{dx^{k}}\phi _{\varepsilon ,\lambda }\left( x\right) =-\frac{%
\beta _{1}\delta ^{-1}}{\sqrt{\lambda }}\dfrac{d^{k}}{dx^{k}}\left( \sin 
\sqrt{\lambda }\left( x-a\right) \right) +O\left( \left( \sqrt{\lambda }%
\right) ^{k-2}e^{\left\vert t\right\vert \left( x-a\right) }\right) , 
\tag{2.23}  \label{2.23}
\end{equation}%
\begin{equation}
\dfrac{d^{k}}{dx^{k}}\phi _{+\varepsilon ,\lambda }\left( x\right) =-\frac{%
\beta _{1}\gamma ^{-1}}{\sqrt{\lambda }}\dfrac{d^{k}}{dx^{k}}\left( \sin 
\sqrt{\lambda }\left( x-a\right) \right) +O\left( \left( \sqrt{\lambda }%
\right) ^{k-2}e^{\left\vert t\right\vert \left( x-a\right) }\right) , 
\tag{2.24}  \label{2.24}
\end{equation}%
if $\beta _{2}=0.$

Then we obtain four distinct cases for the asymptotic behaviour of $\omega
\left( \lambda \right) $ as $\left\vert \lambda \right\vert \rightarrow
\infty ,$ namely;%
\begin{equation}
\omega \left( \lambda \right) =\left\{ 
\begin{array}{l}
\lambda \sqrt{\lambda }\alpha _{1}^{\prime }\beta _{2}\gamma \sin \sqrt{%
\lambda }\left( b-a\right) +O\left( \lambda e^{\left\vert t\right\vert
\left( b-a\right) }\right) ,\text{ \ \ if }\beta _{2}\neq 0,\alpha
_{1}^{\prime }\neq 0, \\ 
\\ 
\lambda \alpha _{2}^{\prime }\beta _{2}\gamma \cos \sqrt{\lambda }\left(
b-a\right) +O\left( \sqrt{\lambda }e^{\left\vert t\right\vert \left(
b-a\right) }\right) ,\text{ \ \ \ if }\beta _{2}\neq 0,\alpha _{1}^{\prime
}=0, \\ 
\\ 
\lambda \alpha _{1}^{\prime }\beta _{1}\gamma \cos \sqrt{\lambda }\left(
b-a\right) +O\left( \sqrt{\lambda }e^{\left\vert t\right\vert \left(
b-a\right) }\right) ,\text{ \ \ \ if }\beta _{2}=0,\alpha _{1}^{\prime }\neq
0, \\ 
\\ 
-\sqrt{\lambda }\alpha _{2}^{\prime }\beta _{1}\gamma \sin \sqrt{\lambda }%
\left( b-a\right) +O\left( e^{\left\vert t\right\vert \left( b-a\right)
}\right) ,\text{ \ \ \ \ if }\beta _{2}=0,\alpha _{1}^{\prime }=0.%
\end{array}%
\right.  \tag{2.25}  \label{2.25}
\end{equation}

Consequently if $\lambda _{0}<\lambda _{1}<...,$ are the zeros of $\omega
\left( \lambda \right) ,$ then we have for sufficiently large $n$ the
following asymptotic formulas%
\begin{equation}
\sqrt{\lambda _{n}}=\left\{ 
\begin{array}{l}
\frac{\left( n-1\right) \pi }{b-a}+O\left( n^{-1}\right) ,\text{ \ \ \ \ \ \
if }\beta _{2}\neq 0,\alpha _{1}^{\prime }\neq 0,\text{\ \ } \\ 
\\ 
\frac{\left( n-1/2\right) \pi }{b-a}+O\left( n^{-1}\right) ,\text{ \ \ \ \
if }\beta _{2}\neq 0,\alpha _{1}^{\prime }=0, \\ 
\\ 
\frac{\left( n-1/2\right) \pi }{b-a}+O\left( n^{-1}\right) ,\text{ \ \ \ \
if }\beta _{2}=0,\alpha _{1}^{\prime }\neq 0, \\ 
\\ 
\frac{n\pi }{b-a}+O\left( n^{-1}\right) ,\text{ \ \ \ \ \ \ \ \ \ \ if }%
\beta _{2}=0,\alpha _{1}^{\prime }=0.%
\end{array}%
\right.  \tag{2.26}  \label{2.26}
\end{equation}

\section{Green Function}

To study the completeness of the eigenvectors of $A,$ and hence the
completeness of the eigenfunctions of the problem (1.1)-(1.7), we construct
the resolvent of $A$ as well as Green's function of the problem (1.1)-(1.7).
We assume without any loss of generality that $\lambda =0$ is not an
eigenvalue of $A.$ Now let $\lambda \in 
\mathbb{C}
$ not be an eigenvalue of $A$ and consider the inhomogenous problem for f$%
\left( x\right) =\left( 
\begin{array}{c}
f\left( x\right) \\ 
f_{1}%
\end{array}%
\right) \in H,$ u$\left( x\right) =\left( 
\begin{array}{c}
u\left( x\right) \\ 
R^{\prime }\left( u\right)%
\end{array}%
\right) \in D\left( A\right) ,$%
\begin{equation}
\left( \lambda \text{I}-A\right) \text{u}\left( x\right) =\text{f}\left(
x\right) ,\text{ }x\in I,  \tag{3.1}  \label{3.1}
\end{equation}%
and I is the identity operator. Since%
\begin{equation}
\left( \lambda \text{I}-A\right) \text{u}\left( x\right) =\lambda \left( 
\begin{array}{c}
u\left( x\right) \\ 
R^{\prime }\left( u\right)%
\end{array}%
\right) -\left( 
\begin{array}{c}
\tau \left( u\right) \\ 
-R\left( u\right)%
\end{array}%
\right) =\left( 
\begin{array}{c}
f\left( x\right) \\ 
f_{1}%
\end{array}%
\right)  \tag{3.2}  \label{3.2}
\end{equation}%
then we have%
\begin{equation}
\left( \lambda \text{I}-\tau \right) \text{u}\left( x\right) =\text{f}\left(
x\right) ,\text{ }x\in I,  \tag{3.3}  \label{3.3}
\end{equation}%
\begin{equation}
\lambda R^{\prime }\left( u\right) +R\left( u\right) =f_{1}.  \tag{3.4}
\label{3.4}
\end{equation}

Now we can represent the general solution of homogeneous differential
equation (1.1), appropriate to equation (3.3) in the following form:%
\begin{equation*}
u\left( x,\lambda \right) =\left\{ 
\begin{array}{c}
c_{1}\phi _{-\varepsilon ,\lambda }\left( x\right) +c_{2}\chi _{-\varepsilon
,\lambda }\left( x\right) ,\text{ \ \ }x\in \left[ a,\theta _{-\varepsilon
}\right) , \\ 
\\ 
c_{3}\phi _{\varepsilon ,\lambda }\left( x\right) +c_{4}\chi _{\varepsilon
,\lambda }\left( x\right) ,\text{ \ \ }x\in \left( \theta _{-\varepsilon
},\theta _{+\varepsilon }\right) , \\ 
\\ 
c_{5}\phi _{+\varepsilon ,\lambda }\left( x\right) +c_{6}\chi _{+\varepsilon
,\lambda }\left( x\right) ,\text{ \ \ \ }x\in \left( \theta _{+\varepsilon
},b\right] ,%
\end{array}%
\right.
\end{equation*}%
in which $c_{i}$ $\left( i=\overline{1,6}\right) $ are arbitrary constants.
By applying the method of variation of the constants, we shall search the
general solution of the non-homogeneous linear differential equation (3.3)
in the following form:%
\begin{equation}
u\left( x,\lambda \right) =\left\{ 
\begin{array}{c}
c_{1}\left( x,\lambda \right) \phi _{-\varepsilon ,\lambda }\left( x\right)
+c_{2}\left( x,\lambda \right) \chi _{-\varepsilon ,\lambda }\left( x\right)
,\text{ \ \ }x\in \left[ a,\theta _{-\varepsilon }\right) , \\ 
\\ 
c_{3}\left( x,\lambda \right) \phi _{\varepsilon ,\lambda }\left( x\right)
+c_{4}\left( x,\lambda \right) \chi _{\varepsilon ,\lambda }\left( x\right) ,%
\text{ \ \ }x\in \left( \theta _{-\varepsilon },\theta _{+\varepsilon
}\right) , \\ 
\\ 
c_{5}\left( x,\lambda \right) \phi _{+\varepsilon ,\lambda }\left( x\right)
+c_{6}\left( x,\lambda \right) \chi _{+\varepsilon ,\lambda }\left( x\right)
,\text{ \ \ \ }x\in \left( \theta _{+\varepsilon },b\right] ,%
\end{array}%
\right.  \tag{3.5}  \label{3.5}
\end{equation}%
where the functions $c_{i}\left( x,\lambda \right) $ $\left( i=\overline{1,6}%
\right) $ satisfy the linear system of equation%
\begin{equation}
\left\{ 
\begin{array}{l}
c_{1}^{\prime }\left( x,\lambda \right) \phi _{-\varepsilon ,\lambda }\left(
x\right) +c_{2}^{\prime }\left( x,\lambda \right) \chi _{-\varepsilon
,\lambda }\left( x\right) =0, \\ 
\\ 
c_{1}^{\prime }\left( x,\lambda \right) \phi _{-\varepsilon ,\lambda
}^{\prime }\left( x\right) +c_{2}^{\prime }\left( x,\lambda \right) \chi
_{-\varepsilon ,\lambda }^{\prime }\left( x\right) =f\left( x\right) ,%
\end{array}%
\right. \text{ \ for \ }x\in \left[ a,\theta _{-\varepsilon }\right) 
\tag{3.6}  \label{3.6}
\end{equation}%
\begin{equation}
\left\{ 
\begin{array}{l}
c_{3}^{\prime }\left( x,\lambda \right) \phi _{\varepsilon ,\lambda }\left(
x\right) +c_{4}^{\prime }\left( x,\lambda \right) \chi _{\varepsilon
,\lambda }\left( x\right) =0, \\ 
\\ 
c_{3}^{\prime }\left( x,\lambda \right) \phi _{\varepsilon ,\lambda
}^{\prime }\left( x\right) +c_{4}^{\prime }\left( x,\lambda \right) \chi
_{\varepsilon ,\lambda }^{\prime }\left( x\right) =f\left( x\right) ,%
\end{array}%
\right. \text{ \ for \ }x\in \left( \theta _{-\varepsilon },\theta
_{+\varepsilon }\right)  \tag{3.7}  \label{3.7}
\end{equation}%
\begin{equation}
\left\{ 
\begin{array}{l}
c_{5}^{\prime }\left( x,\lambda \right) \phi _{+\varepsilon ,\lambda }\left(
x\right) +c_{6}^{\prime }\left( x,\lambda \right) \chi _{+\varepsilon
,\lambda }\left( x\right) =0, \\ 
\\ 
c_{5}^{\prime }\left( x,\lambda \right) \phi _{+\varepsilon ,\lambda
}^{\prime }\left( x\right) +c_{6}^{\prime }\left( x,\lambda \right) \chi
_{+\varepsilon ,\lambda }^{\prime }\left( x\right) =f\left( x\right) ,%
\end{array}%
\right. \text{ \ for \ }x\in \left( \theta _{+\varepsilon },b\right] . 
\tag{3.8}  \label{3.8}
\end{equation}

Since $\lambda $ is not an eigenvalue and $\omega _{-\varepsilon }\left(
\lambda \right) \neq 0,$ $\omega _{\varepsilon }\left( \lambda \right) \neq
0,$ $\omega _{+\varepsilon }\left( \lambda \right) \neq 0,$ each of the
linear systems in (3.6)-(3.8) have a unique solution which leads%
\begin{equation}
\left\{ 
\begin{array}{l}
c_{1}\left( x,\lambda \right) =\frac{1}{\omega _{-\varepsilon }\left(
\lambda \right) }\underset{x}{\overset{\theta _{-\varepsilon }}{\int }}\chi
_{-\varepsilon ,\lambda }\left( y\right) f\left( y\right) dy+c_{1}\left(
\lambda \right) , \\ 
c_{2}\left( x,\lambda \right) =\frac{1}{\omega _{-\varepsilon }\left(
\lambda \right) }\underset{a}{\overset{x}{\int }}\phi _{-\varepsilon
,\lambda }\left( y\right) f\left( y\right) dy+c_{2}\left( \lambda \right) ,%
\end{array}%
\right. \text{ for \ }x\in \left[ a,\theta _{-\varepsilon }\right)  \tag{3.9}
\label{3.9}
\end{equation}%
\begin{equation}
\left\{ 
\begin{array}{l}
c_{3}\left( x,\lambda \right) =\frac{1}{\omega _{\varepsilon }\left( \lambda
\right) }\underset{x}{\overset{\theta _{+\varepsilon }}{\int }}\chi
_{\varepsilon ,\lambda }\left( y\right) f\left( y\right) dy+c_{3}\left(
\lambda \right) , \\ 
c_{4}\left( x,\lambda \right) =\frac{1}{\omega _{\varepsilon }\left( \lambda
\right) }\underset{\theta _{-\varepsilon }}{\overset{x}{\int }}\phi
_{\varepsilon ,\lambda }\left( y\right) f\left( y\right) dy+c_{4}\left(
\lambda \right) ,%
\end{array}%
\right. \text{ \ for \ }x\in \left( \theta _{-\varepsilon },\theta
_{+\varepsilon }\right)  \tag{3.10}  \label{3.10}
\end{equation}%
\begin{equation}
\left\{ 
\begin{array}{l}
c_{5}\left( x,\lambda \right) =\frac{1}{\omega _{+\varepsilon }\left(
\lambda \right) }\underset{x}{\overset{b}{\int }}\chi _{+\varepsilon
,\lambda }\left( y\right) f\left( y\right) dy+c_{5}\left( \lambda \right) ,
\\ 
c_{6}\left( x,\lambda \right) =\frac{1}{\omega _{+\varepsilon }\left(
\lambda \right) }\underset{\theta _{+\varepsilon }}{\overset{x}{\int }}\phi
_{+\varepsilon ,\lambda }\left( y\right) f\left( y\right) dy+c_{6}\left(
\lambda \right) ,%
\end{array}%
\right. \text{ \ for \ }x\in \left( \theta _{+\varepsilon },b\right] 
\tag{3.11}  \label{3.11}
\end{equation}%
where $c_{i}\left( \lambda \right) $ $\left( i=\overline{1,6}\right) $ are
arbitrary constants. Substituting (3.9)-(3.11) into (3.5), we obtain the
solution of the equation (3.3)%
\begin{equation}
u\left( x,\lambda \right) =\left\{ 
\begin{array}{l}
\frac{\phi _{-\varepsilon ,\lambda }\left( x\right) }{\omega _{-\varepsilon
}\left( \lambda \right) }\underset{x}{\overset{\theta _{-\varepsilon }}{\int 
}}\chi _{-\varepsilon ,\lambda }\left( y\right) f\left( y\right) dy+\frac{%
\chi _{-\varepsilon ,\lambda }\left( x\right) }{\omega _{-\varepsilon
}\left( \lambda \right) }\underset{a}{\overset{x}{\int }}\phi _{-\varepsilon
,\lambda }\left( y\right) f\left( y\right) dy+ \\ 
c_{1}\left( \lambda \right) \phi _{-\varepsilon ,\lambda }\left( x\right)
+c_{2}\left( \lambda \right) \chi _{-\varepsilon ,\lambda }\left( x\right) ,%
\text{ \ \ \ \ \ \ \ \ \ \ \ \ \ \ \ \ \ \ \ \ }x\in \left[ a,\theta
_{-\varepsilon }\right) , \\ 
\\ 
\frac{\phi _{\varepsilon ,\lambda }\left( x\right) }{\omega _{\varepsilon
}\left( \lambda \right) }\underset{x}{\overset{\theta _{+\varepsilon }}{\int 
}}\chi _{\varepsilon ,\lambda }\left( y\right) f\left( y\right) dy+\frac{%
\chi _{\varepsilon ,\lambda }\left( x\right) }{\omega _{\varepsilon }\left(
\lambda \right) }\underset{\theta _{-\varepsilon }}{\overset{x}{\int }}\phi
_{\varepsilon ,\lambda }\left( y\right) f\left( y\right) dy+ \\ 
c_{3}\left( \lambda \right) \phi _{\varepsilon ,\lambda }\left( x\right)
+c_{4}\left( \lambda \right) \chi _{\varepsilon ,\lambda }\left( x\right) ,%
\text{ \ \ \ \ \ \ \ \ \ \ \ \ \ \ \ \ \ \ \ \ }x\in \left( \theta
_{-\varepsilon },\theta _{+\varepsilon }\right) , \\ 
\\ 
\frac{\phi _{+\varepsilon ,\lambda }\left( x\right) }{\omega _{+\varepsilon
}\left( \lambda \right) }\underset{x}{\overset{b}{\int }}\chi _{+\varepsilon
,\lambda }\left( y\right) f\left( y\right) dy+\frac{\chi _{+\varepsilon
,\lambda }\left( x\right) }{\omega _{+\varepsilon }\left( \lambda \right) }%
\underset{\theta _{+\varepsilon }}{\overset{x}{\int }}\phi _{+\varepsilon
,\lambda }\left( y\right) f\left( y\right) dy+ \\ 
c_{5}\left( \lambda \right) \phi _{+\varepsilon ,\lambda }\left( x\right)
+c_{6}\left( \lambda \right) \chi _{+\varepsilon ,\lambda }\left( x\right) ,%
\text{ \ \ \ \ \ \ \ \ \ \ \ \ \ \ \ \ \ \ \ \ }x\in \left( \theta
_{+\varepsilon },b\right] .%
\end{array}%
\right.  \tag{3.12}  \label{3.12}
\end{equation}

Then from the boundary conditions (3.4) and (1.2) and the transmission
conditions (1.4)-(1.7) we get%
\begin{equation}
\begin{array}{l}
c_{1}\left( \lambda \right) =\frac{1}{\omega _{\varepsilon }\left( \lambda
\right) }\underset{\theta _{-\varepsilon }}{\overset{\theta _{+\varepsilon }}%
{\int }}\chi _{\varepsilon ,\lambda }\left( y\right) f\left( y\right) dy+%
\frac{1}{\omega _{+\varepsilon }\left( \lambda \right) }\underset{\theta
_{+\varepsilon }}{\overset{b}{\int }}\chi _{+\varepsilon ,\lambda }\left(
y\right) f\left( y\right) dy+\frac{f_{1}}{\omega _{+\varepsilon }\left(
\lambda \right) }, \\ 
\\ 
c_{2}\left( \lambda \right) =0,\text{ \ \ \ \ \ \ }c_{3}\left( \lambda
\right) =\frac{1}{\omega _{+\varepsilon }\left( \lambda \right) }\underset{%
\theta _{+\varepsilon }}{\overset{b}{\int }}\chi _{+\varepsilon ,\lambda
}\left( y\right) f\left( y\right) dy+\frac{f_{1}}{\omega _{+\varepsilon
}\left( \lambda \right) }, \\ 
\\ 
c_{4}\left( \lambda \right) =\frac{1}{\omega _{-\varepsilon }\left( \lambda
\right) }\underset{a}{\overset{\theta _{-\varepsilon }}{\int }}\phi
_{-\varepsilon ,\lambda }\left( y\right) f\left( y\right) dy,\text{ \ \ \ \
\ \ \ \ \ \ }c_{5}\left( \lambda \right) =\frac{f_{1}}{\omega _{+\varepsilon
}\left( \lambda \right) }, \\ 
\\ 
c_{6}\left( \lambda \right) =\frac{1}{\omega _{-\varepsilon }\left( \lambda
\right) }\underset{a}{\overset{\theta _{-\varepsilon }}{\int }}\phi
_{-\varepsilon ,\lambda }\left( y\right) f\left( y\right) dy+\frac{1}{\omega
_{\varepsilon }\left( \lambda \right) }\underset{\theta _{-\varepsilon }}{%
\overset{\theta _{+\varepsilon }}{\int }}\phi _{\varepsilon ,\lambda }\left(
y\right) f\left( y\right) dy.%
\end{array}
\tag{3.13}  \label{3.13}
\end{equation}

Substituting (3.13) and (2.11) into (3.12), then (3.12) can be written as%
\begin{equation}
u\left( x,\lambda \right) =\left\{ 
\begin{array}{l}
\frac{\phi _{-\varepsilon ,\lambda }\left( x\right) }{\omega \left( \lambda
\right) }\underset{x}{\overset{\theta _{-\varepsilon }}{\int }}\chi
_{-\varepsilon ,\lambda }\left( y\right) f\left( y\right) dy+\frac{\chi
_{-\varepsilon ,\lambda }\left( x\right) }{\omega \left( \lambda \right) }%
\underset{a}{\overset{x}{\int }}\phi _{-\varepsilon ,\lambda }\left(
y\right) f\left( y\right) dy+ \\ 
\frac{\delta ^{2}\phi _{-\varepsilon ,\lambda }\left( x\right) }{\omega
\left( \lambda \right) }\underset{\theta _{-\varepsilon }}{\overset{\theta
_{+\varepsilon }}{\int }}\chi _{\varepsilon ,\lambda }\left( y\right)
f\left( y\right) dy+\frac{\delta ^{2}\phi _{-\varepsilon ,\lambda }\left(
x\right) }{\omega \left( \lambda \right) }\underset{\theta _{+\varepsilon }}{%
\overset{b}{\int }}\chi _{+\varepsilon ,\lambda }\left( y\right) f\left(
y\right) dy+ \\ 
\frac{\gamma ^{2}f_{1}}{\omega \left( \lambda \right) }\phi _{-\varepsilon
,\lambda }\left( x\right) ,\text{\ \ \ \ \ \ \ \ \ \ \ \ \ \ \ \ \ \ \ \ \ \
\ \ \ \ \ \ \ \ \ \ \ \ \ \ \ \ \ \ \ \ \ \ \ \ \ \ \ \ \ \ }x\in \left[
a,\theta _{-\varepsilon }\right) , \\ 
\\ 
\frac{\delta ^{2}\phi _{\varepsilon ,\lambda }\left( x\right) }{\omega
\left( \lambda \right) }\underset{x}{\overset{\theta _{+\varepsilon }}{\int }%
}\chi _{\varepsilon ,\lambda }\left( y\right) f\left( y\right) dy+\frac{%
\delta ^{2}\chi _{\varepsilon ,\lambda }\left( x\right) }{\omega \left(
\lambda \right) }\underset{\theta _{-\varepsilon }}{\overset{x}{\int }}\phi
_{\varepsilon ,\lambda }\left( y\right) f\left( y\right) dy+ \\ 
\frac{\gamma ^{2}\phi _{\varepsilon ,\lambda }\left( x\right) }{\omega
\left( \lambda \right) }\underset{\theta _{+\varepsilon }}{\overset{b}{\int }%
}\chi _{+\varepsilon ,\lambda }\left( y\right) f\left( y\right) dy+\frac{%
\chi _{\varepsilon ,\lambda }\left( x\right) }{\omega \left( \lambda \right) 
}\underset{a}{\overset{\theta _{-\varepsilon }}{\int }}\phi _{-\varepsilon
,\lambda }\left( y\right) f\left( y\right) dy+ \\ 
\frac{\gamma ^{2}f_{1}}{\omega \left( \lambda \right) }\phi _{\varepsilon
,\lambda }\left( x\right) ,\text{ \ \ \ \ \ \ \ \ \ \ \ \ \ \ \ \ \ \ \ \ \
\ \ \ \ \ \ \ \ \ \ \ \ \ \ \ \ \ \ \ \ \ \ \ \ \ \ \ \ \ }x\in \left(
\theta _{-\varepsilon },\theta _{+\varepsilon }\right) , \\ 
\\ 
\frac{\gamma ^{2}\phi _{+\varepsilon ,\lambda }\left( x\right) }{\omega
\left( \lambda \right) }\underset{x}{\overset{b}{\int }}\chi _{+\varepsilon
,\lambda }\left( y\right) f\left( y\right) dy+\frac{\gamma ^{2}\chi
_{+\varepsilon ,\lambda }\left( x\right) }{\omega \left( \lambda \right) }%
\underset{\theta _{+\varepsilon }}{\overset{x}{\int }}\phi _{+\varepsilon
,\lambda }\left( y\right) f\left( y\right) dy+ \\ 
\frac{\chi _{+\varepsilon ,\lambda }\left( x\right) }{\omega \left( \lambda
\right) }\underset{a}{\overset{\theta _{-\varepsilon }}{\int }}\phi
_{-\varepsilon ,\lambda }\left( y\right) f\left( y\right) dy+\frac{\delta
^{2}\chi _{+\varepsilon ,\lambda }\left( x\right) }{\omega \left( \lambda
\right) }\underset{\theta _{-\varepsilon }}{\overset{\theta _{+\varepsilon }}%
{\int }}\phi _{\varepsilon ,\lambda }\left( y\right) f\left( y\right) dy+ \\ 
\frac{\gamma ^{2}f_{1}}{\omega \left( \lambda \right) }\phi _{+\varepsilon
,\lambda }\left( x\right) ,\text{ \ \ \ \ \ \ \ \ \ \ \ \ \ \ \ \ \ \ \ \ \
\ \ \ \ \ \ \ \ \ \ \ \ \ \ \ \ \ \ \ \ \ \ \ \ \ \ \ \ \ }x\in \left(
\theta _{+\varepsilon },b\right] .%
\end{array}%
\right.  \tag{3.14}  \label{3.14}
\end{equation}

Hence, we have%
\begin{eqnarray}
\text{u}\left( x\right) &=&\left( \lambda \text{I}-A\right) ^{-1}\text{f}%
\left( x\right)  \notag \\
&=&\left( \underset{{\Large R}^{\prime }\left( u\right) }{%
\begin{array}{c}
\underset{a}{\overset{\theta _{-\varepsilon }}{\int }}G\left( x,y,\lambda
\right) f\left( y\right) dy+\delta ^{2}\underset{\theta _{-\varepsilon }}{%
\overset{\theta _{+\varepsilon }}{\int }}G\left( x,y,\lambda \right) f\left(
y\right) dy+ \\ 
\gamma ^{2}\underset{\theta _{+\varepsilon }}{\overset{b}{\int }}G\left(
x,y,\lambda \right) f\left( y\right) dy+\dfrac{\gamma ^{2}}{\omega \left(
\lambda \right) }f_{1}\phi _{\lambda }\left( x\right) \text{ \ \ \ \ \ \ \ \
\ \ \ \ \ }%
\end{array}%
}\right)  \TCItag{3.15}  \label{3.15}
\end{eqnarray}%
where%
\begin{equation}
G\left( x,y,\lambda \right) =\frac{1}{\omega \left( \lambda \right) }\left\{ 
\begin{array}{l}
\chi _{\lambda }\left( x\right) \phi _{-\varepsilon ,\lambda }\left(
y\right) ,\text{ \ \ \ \ \ }a\leq y\leq x\leq \theta _{-\varepsilon }, \\ 
\\ 
\chi _{\lambda }\left( y\right) \phi _{-\varepsilon ,\lambda }\left(
x\right) ,\text{ \ \ \ \ \ }a\leq x\leq y\leq \theta _{-\varepsilon }, \\ 
\\ 
\chi _{\lambda }\left( x\right) \phi _{\lambda }\left( y\right) ,\text{ \ \
\ \ \ }\theta _{-\varepsilon }\leq y\leq x\leq \theta _{+\varepsilon }, \\ 
\\ 
\chi _{\lambda }\left( y\right) \phi _{\lambda }\left( x\right) ,\text{ \ \
\ \ \ }\theta _{-\varepsilon }\leq x\leq y\leq \theta _{+\varepsilon }, \\ 
\\ 
\chi _{+\varepsilon ,\lambda }\left( x\right) \phi _{\lambda }\left(
y\right) ,\text{ \ \ \ \ \ }\theta _{+\varepsilon }\leq y\leq x\leq b, \\ 
\\ 
\chi _{+\varepsilon ,\lambda }\left( y\right) \phi _{\lambda }\left(
x\right) ,\text{ \ \ \ \ \ }\theta _{+\varepsilon }\leq x\leq y\leq b,%
\end{array}%
\right.  \tag{3.16}  \label{3.16}
\end{equation}%
is Green's function of the problem (1.1)-(1.7).

\section{The Sampling Theorem}

In this section we derive two sampling theorems associated with the problem
(1.1)-(1.7). For convenience we may assume that the eigenvectors of $A$ are
real valued.

\begin{theorem}
Consider the problem (1.1)-(1.7), and let%
\begin{equation}
\phi _{\lambda }\left( x\right) =\left\{ 
\begin{array}{l}
\phi _{-\varepsilon ,\lambda }\left( x\right) ,\text{ \ \ }x\in \left[
a,\theta _{-\varepsilon }\right) , \\ 
\phi _{\varepsilon ,\lambda }\left( x\right) ,\text{ \ }x\in \left( \theta
_{-\varepsilon },\theta _{+\varepsilon }\right) , \\ 
\phi _{+\varepsilon ,\lambda }\left( x\right) ,\text{ \ \ }x\in \left(
\theta _{+\varepsilon },b\right] ,%
\end{array}%
\right.  \tag{4.1}  \label{4.1}
\end{equation}%
\textit{be the solution defined above. Let }$g\left( .\right) \in
L_{2}\left( a,b\right) $\textit{\ and }%
\begin{equation}
F\left( \lambda \right) =\underset{a}{\overset{\theta _{-\varepsilon }}{\int 
}}g\left( x\right) \phi _{-\varepsilon ,\lambda }\left( x\right) dx+\delta
^{2}\underset{\theta _{-\varepsilon }}{\overset{\theta _{+\varepsilon }}{%
\int }}g\left( x\right) \phi _{\varepsilon ,\lambda }\left( x\right)
dx+\gamma ^{2}\underset{\theta _{+\varepsilon }}{\overset{b}{\int }}g\left(
x\right) \phi _{+\varepsilon ,\lambda }\left( x\right) dx.  \tag{4.2}
\label{4.2}
\end{equation}
\end{theorem}

\textit{Then }$F\left( \lambda \right) $\textit{\ is an entire function of
exponential type }$(b-a)$\textit{\ that can be reconstructed from its values
at the points }$\left\{ \lambda _{n}\right\} _{n=0}^{\infty }$\textit{\ \ \
via the sampling formula}%
\begin{equation}
F\left( \lambda \right) =\overset{\infty }{\underset{n=0}{\sum }}F\left(
\lambda _{n}\right) \frac{\omega \left( \lambda \right) }{\left( \lambda
-\lambda _{n}\right) \omega ^{\prime }\left( \lambda _{n}\right) }. 
\tag{4.3}  \label{4.3}
\end{equation}

\textit{The series (4.2) converges absolutely on }$%
\mathbb{C}
$\textit{\ and uniformly on compact subset of }$%
\mathbb{C}
.$\textit{\ Here }$\omega \left( \lambda \right) $\textit{\ is the entire
function defined in (2.11).}

\begin{proof}
...
\end{proof}

\begin{remark}
To see that expansion (4.3) is a Lagrange type interpolation, we may replace 
$\omega \left( \lambda \right) $ by the canonical product%
\begin{equation}
\widetilde{\omega }\left( \lambda \right) =\left\{ 
\begin{array}{c}
\overset{\infty }{\underset{n=0}{\Pi }}\left( 1-\dfrac{\lambda }{\lambda _{n}%
}\right) ,\text{ if zero is not an eigenvalue,} \\ 
\\ 
\lambda \overset{\infty }{\underset{n=1}{\Pi }}\left( 1-\dfrac{\lambda }{%
\lambda _{n}}\right) ,\text{ if }\lambda _{0}=0\text{ is an eigenvalue.}%
\end{array}%
\right.  \tag{4.26}  \label{4.26}
\end{equation}
\end{remark}

From Hadamard's factorization theorem, see $\left[ 31\right] ,$ $\omega
\left( \lambda \right) =h\left( \lambda \right) \widetilde{\omega }\left(
\lambda \right) ,$ where $h\left( \lambda \right) $ is an entire function
with no zeros. Thus,%
\begin{equation}
\frac{\omega \left( \lambda \right) }{\omega ^{\prime }\left( \lambda
_{n}\right) }=\frac{h\left( \lambda \right) \widetilde{\omega }\left(
\lambda \right) }{h\left( \lambda _{n}\right) \widetilde{\omega }^{\prime
}\left( \lambda _{n}\right) }  \tag{4.27}  \label{4.27}
\end{equation}%
and (4.2), (4.3) remain valid for the function $F\left( \lambda \right)
/h\left( \lambda \right) .$ Hence%
\begin{equation}
F\left( \lambda \right) =\overset{\infty }{\underset{n=0}{\sum }}F\left(
\lambda _{n}\right) \frac{h\left( \lambda \right) \widetilde{\omega }\left(
\lambda \right) }{\left( \lambda -\lambda _{n}\right) h\left( \lambda
_{n}\right) \widetilde{\omega }^{\prime }\left( \lambda _{n}\right) }. 
\tag{4.28}  \label{4.28}
\end{equation}

We may redefine (4.2) by taking kernel $\phi _{\lambda }\left( .\right)
/h\left( \lambda \right) =\widetilde{\phi }_{\lambda }\left( .\right) $ to
get%
\begin{equation}
\widetilde{F}\left( \lambda \right) =\frac{F\left( \lambda \right) }{h\left(
\lambda \right) }=\overset{\infty }{\underset{n=0}{\sum }}\widetilde{F}%
\left( \lambda _{n}\right) \frac{\widetilde{\omega }\left( \lambda \right) }{%
\left( \lambda -\lambda _{n}\right) \widetilde{\omega }^{\prime }\left(
\lambda _{n}\right) }.  \tag{4.29}  \label{4.29}
\end{equation}

The next theorem is devoted to give interpolation sampling expansions
associated with the problem (1.1)-(1.7) for integral transforms whose
kernels defined in terms of Green's function (see $\left[ 20,32\right] )$.
As we see in (3.16), Green' function $G\left( x,y,\lambda \right) $ of the
problem (1.1)-(1.7) has simple poles at $\left\{ \lambda _{n}\right\}
_{n=0}^{\infty }.$

Define the function $G\left( x,\lambda \right) $ to be $G\left( x,\lambda
\right) :=\omega \left( \lambda \right) G\left( x,y_{0},\lambda \right) ,$
where $y_{0}\in I$ is a fixed point and $\omega \left( \lambda \right) $ is
the function defined in equation (2.11) or it is the canonical product
(4.26).

\begin{theorem}
Let g$\left( .\right) \in L_{2}\left( a,b\right) $ and $F\left( \lambda
\right) $ the integral transform%
\begin{equation}
F\left( \lambda \right) =\underset{a}{\overset{\theta _{-\varepsilon }}{\int 
}}G\left( x,\lambda \right) \overline{g}\left( x\right) dx+\delta ^{2}%
\underset{\theta _{-\varepsilon }}{\overset{\theta _{+\varepsilon }}{\int }}%
G\left( x,\lambda \right) \overline{g}\left( x\right) dx+\gamma ^{2}\underset%
{\theta _{+\varepsilon }}{\overset{b}{\int }}G\left( x,\lambda \right) 
\overline{g}\left( x\right) dx.  \tag{4.30}  \label{4.30}
\end{equation}
\end{theorem}

\textit{Then }$F\left( \lambda \right) $\textit{\ is an entire function of
exponential type }$(b-a)$\textit{\ which admits the sampling representation}%
\begin{equation}
F\left( \lambda \right) =\overset{\infty }{\underset{n=0}{\sum }}F\left(
\lambda _{n}\right) \frac{\omega \left( \lambda \right) }{\left( \lambda
-\lambda _{n}\right) \omega ^{\prime }\left( \lambda _{n}\right) }. 
\tag{4.31}  \label{4.31}
\end{equation}

\textit{Series (4.31) converges absolutely on }$%
\mathbb{C}
$\textit{\ and uniformly on compact subset of }$%
\mathbb{C}
.$\textit{\ }

\begin{proof}
...
\end{proof}

\end{document}